\documentstyle[12pt]{article}
\begin{document}
\title{ Differentiation of Integrals
}
\author{{    Shunchao Long }}
\date{}
\maketitle
\begin{center}
\begin{minipage}{120mm}
\vskip 0.1in {
\begin{center}
{\bf Abstract}
\end{center}
\par
 No functions class
 for general measurable sets classes are known whose functions have  the property of differentiability of integrals   associated to such sets classes.

\par In this paper,
 we give  some subspaces of $L^s$ with $1<s<\infty$, whose functions are proven to have the  differentiability of integrals associated to measurable sets classes in ${\bf R}^n $,  this gives an answer to a question stated by Stein in his book   Harmonic Analysis.
  We give also a example of  some functions   in these classes on ${\bf R}^2 $, which is continuous nowhere.

}
\end{minipage}
\end{center}
\vskip 0.1in
\baselineskip 16pt
\par

 {\bf  1.  Introduction }
 It is considered that the question of the differentiability of integrals in ${\bf R}^n $
represents one of the main issues in real-variable theory, Stein  formulated this question in\cite{S3} as follows:

  \par {\bf Problem A}
for what collections $\mathcal{C}$ of sets $\{ C\}$,
is it true that for "all" $f$
$$
\lim_{{\rm diam}(C)\rightarrow 0,
\\x\in C\in \mathcal{C}} \frac{1}{|C|}\int_C f(y)dy =f(x) ~~{\rm a.e.} ~~x\in {\bf R}^n ?\eqno (1)
$$
 \par {\bf Problem B} given $\mathcal{C}$, for what classes of functions $f$ does   (1) hold?

\par  If     $   {\mathcal{C}}    $ $ =\{ all ~ cubes\}$ or $\{ all ~ balls\}$, the   Lebesgue's differentiation theorem asserts that (1) holds for all  $f$ in $L^1_{{\rm loc}}$, (see \cite{ S3}).

\par   If    $   {\mathcal{C}}    $$ =\{rectangles: whose~major~ axes ~ point ~in ~a ~ fixed ~directions\}$,  then (1) holds for all   $f$ in $L^p_{{\rm loc}}$ with $ p>1$, but fails when $p=1$, (see \cite{S3}).

\par  If    $    {\mathcal{C}}    $$ =\{all ~rectangles \}$,   for each $p$ with $1\leq p<\infty $,
there exists an $f \in L^p$ such that (1) does not hold,
 (see \cite{S3}).

\par   If    $   {\mathcal{C}}    $$ =\{tubes~ that~ point~ in~ the~ given~ directions \}$,  then (1) holds for all   $f$ in $L^p_{{\rm loc}}$ with $ p>1$, (see \cite{PR}).

\par   If the set in  $ \mathcal{C}$ have bounded eccentricity, in the sense that the ratio between the smallest ball containing $C$ to the largest ball contained in $C$  is uniformly bounded over $C\in \mathcal{C}$,  then (1) holds for all  $f$ in $L^1_{{\rm loc}}$,
 (see \cite{S3}).

Efforts on  the issues above includes  the works in  a lot of  literatures, for example, \cite{SJ,CF,SJ1,NSW,B}, etc.

  Even if for the class of arbitrary rectangles,
  as we see,  $L^p$ is not suitable for (1). So far, only the   class $C_c$ of continuous  functions with compact supports  for the class of arbitrary rectangles, while no class
 for general measurable sets classes,    have been seen in the literatures, whose  functions $f$ satisfying  (1).

\par Let    $    {\mathcal{C}} _0   =\{$all measurable sets $C\subset {\bf R}^n  $ with $|C|\neq 0\}$.
  In this paper,   we give some functions classes   on ${\bf R}^n $,
  which are subspaces of $L^p$ with $ 1<p<\infty,$  whose functions are proven to satisfy (1) for   ${\mathcal{C}}={\mathcal{C}}_0$. This gives an answer for     Problem B  above.
We give also a example of  some functions in these classes on ${\bf R}^2 $, which is continuous nowhere.

\par
 Let  $    {\mathcal{C}}_1    $$ =\{ cubes  ~ Q
    :|Q|\geq 1   \}$.
Let $  0< p < \infty , - n < \alpha < \infty  $. Denote that
$  \bar{p} = {\rm min} \{p, 1\}
 $.

 \par
 \par {\bf Definition  1}
    (A)    $a(x)$ is said to be a
 $  (p, \alpha) $-block, if
\par~~~~ (i)~~ supp $a\subseteq Q \in{\mathcal{C}} _{1 }, $
 ~~~~ (ii)~~ $\|a\|_{L^{\infty}   }\leq |Q|^{-\alpha/pn-1/p  }.$
\par (B)  A $  (p, \alpha) $-block $a(x)$ is said to be a
 $  (p, \alpha) $-continuous block, if $a(x)\in C_c$.

 \par (C)    $a(x)$ is said to be a
 $  (p, \alpha) $-characteristic block, if $a(x) =|Q|^{-\alpha/pn-1/p} \chi_Q(x) $ with  $  Q \in{\mathcal{C}} _{1 } $

\par {\bf Definition  2}
    (A)
$ BL^{p,   \alpha}
  =\{  f :
 f=\sum _{k=-\infty}^{\infty} \lambda _ka_k ,
 \textrm { where
 each $a_k$ is a $ _{ }(p,  \alpha) $ }$-block, $
 \sum _{k=-\infty}^{\infty} |\lambda _k|^{\bar{p}} <+ \infty \},
  $
 here,  the
 "convergence"  means  a.e. convergence.
 Moreover,  define a quasinorm on $  _{ }BL^{p,   \alpha}
 $ by
$\|f\|_{ _{ }BL^{p,   \alpha}
}=
 \inf \left(\sum _{k=-\infty}^{\infty}|\lambda _k|^{\bar{p}}\right)^{1/{\bar{p}}},$
where the infimum is taken over all the a.e-equal forms of $f=\sum _{k=-\infty}^{\infty} \lambda _ka_k $.

\par (B)
$
  _{ }BC^{p,   \alpha}$ is defined  if   each  $a_k$ in the definition above   is a $(p,  \alpha) $-continuous block.

\par (C)
$
  _{ }B^{p,   \alpha}$ is defined  if  each  $a_k$ in the definition above   is a $(p,  \alpha) $-characteristic block.



 {\bf  2.   The Maximal Theorem  } Let
 $$
M _{ {\mathcal{C}} _0 }f(x)=\sup_{x\in C \in { {\mathcal{C}} _0 }, {\rm diam}(C)<1 } \frac{1}{|C|}\int_C |f(y)|dy. \eqno (2)
$$

 \par {\bf  Theorem 1} Let $0<p<\infty$ and $  - n < \alpha < n(p-1)$.  $M_{ {\mathcal{C}} _0 }$ is bounded on
$
BL^{p,   \alpha} $.

\par We will prove Theorem 1   with the aid of  the idea of
 \cite{ CW}.

  \par {\bf Definition 3}   Let $  0< p <+\infty,- n < \alpha < +\infty
 , -\infty <A_0,B_0 <\infty, A_0-B_0
 = -\frac{1}{p}-\frac{\alpha}{np}$ and $
  -\infty<\varepsilon < A_0. $ Set $
   a=A_0 -\varepsilon ,$ and $b=B_0-\varepsilon $.
  A function $M(x )\in L^{\infty} $
 is said to be a  $ (p,   \alpha, \varepsilon )$-molecular (centered at  $x_0 $), if
\par (i) $M(x )|x-x_0 | ^{nb}\in L^{\infty} ,$ ~~
 (ii)  $\|M\|^{a/b}_{L^{\infty}   }\|M(x )|x-x_0 |^{nb}\|^{1-a/b}_{L^{\infty} }
\equiv \Re (M)<\infty.$



\par {\bf Theorem 2  }    Let $  p, \alpha  ,  \varepsilon  , a ,b$
as in Definition~3.
 Let $M(x)$ be a $ (p,  \alpha, \varepsilon
)$-molecular   centered at any point $x_0\in{\bf R}^{n}$.
  If there exists $Q_{0} \in   {\mathcal{C}}    $$_{1 }$ such that
$$\|M\|^{ }_{L^{\infty}   }\leq
|Q_{0} |^{ -1/p-\alpha/np}, $$
     then
      $M \in BL^{p,   \alpha}$ and
$\|M \|^{ }_{_{ }BL^{p,   \alpha}}\leq C
\Re(M), $
where
the constant $C$ is independent of $M$.

\par {\bf Proof}
 Without loss of generality, we can assume $\Re(M)=1$. In fact,
assume $\|M \|^{ }_{_{ }{B}L^{p,   \alpha}}\leq
C  $ holds whenever $\Re(M)=1$. Then, for general $M,$ let
$M'=M/\Re(M)$. We have  $\Re(M')=1$ and hence $\|\Re(M)M' \|^{
}_{_{ }BL^{p,   \alpha}}\leq \Re(M)\| M' \|^{
}_{_{ }BL^{p,   \alpha}}\leq C \Re(M)$. And we
can also assume that  $M(x)$  has the center at $x_0=0$ by a translation
transformation.

\par
 Define cube $Q_{1} $ centered at
$ 0$ whose sides   parallel
to the axes by setting
$$\|M\|^{ }_{L^{\infty}   }=
|Q_{1}|^{-1/p-\alpha/np}.$$
We see that $|Q_1|\geq |Q_0|$ since $
|Q_{1}|^{-1/p-\alpha/np}\leq |Q_{0}|^{-1/p-\alpha/np}$ and $\alpha>-n$,
so we can assume that $Q_{1} \in   {\mathcal{C}}    $$_{1}$.
Let
  $$ E_0=Q_{1} ,
  E_k=2^{k}Q_{1}\setminus 2^{k-1}Q_{1} ,   k= 1, 2, \cdots         . $$
 Set $M_k=M{\chi _{E_k}}$.
   By $\Re (M)=1$, we have
 \begin{eqnarray*}
 \left\|M(x)
 |x | ^{nb}\right\|_{L^{\infty}
}=\|M\|^{-\frac{a}{b}\frac{b}{b-a}}_{L^{\infty}}
=|Q_{1} |^{-\left( -\frac{1}{p}
-\frac{\alpha}{np} \right)
 \frac{a}{b-a}}
=|Q_{1} |^a.
\end{eqnarray*}
 For  $x\in E_k $, $k=1,2,\cdots $, we have
\begin{eqnarray*}
|M_k(x)|=|M_k(x)|  |x | ^{nb}
 |x | ^{-nb}
 \leq    |2^{k-1}Q_{1}|^{-b} |M_k(x)| |x | ^{nb}
 =  2^{n  b } |2^{k}Q_{1}|^{-b} |M_k(x)| |x | ^{nb},
\end{eqnarray*}
since $b>0$. It follows that
\begin{eqnarray*}
\|M_k\|_{L^{\infty}   }
&\leq & 2^{n  b }|2^{k }Q_{1}|^{-b} \left\|M_k(x)
 |x | ^{nb}\right\|_{L^{\infty}}
 \leq
2^{n  b }|2^{k }Q_{1}|^{-b}|Q_{1} |^a
\\&= &
2^{n  b } 2^{-kna}|2^{k }Q_{1}|^{ -1/p-\alpha/np},
\end{eqnarray*}
for
$k= 1,2,...$ , while
$$\|M_0\|^{ }_{L^{\infty}   }\leq  \|M\|^{ }_{L^{\infty}   }= |Q_{1} |^{ -1/p-\alpha/np} \leq 2^{n  b }|Q_{1} |^{ -1/p-\alpha/np}
.$$
Let
$$
a_k(x)=2^{-n  b }2^{kna}M_k(x), ~~k=0,1,2,...,
$$
 then
$$
M(x)=\sum \limits_{k=0}^{\infty}M_k(x)=\sum \limits_{k=0}^{\infty}
2^{n  b } 2^{-kna}2^{-n  b }2^{kna}
 M_k(x)=2^{n  b } \sum \limits_{k=0}^{\infty}2^{-kna}a_k(x),$$
and each $a_k$ is a $ _{ }(p, \alpha)$-block centered at  $ 0$ with
${\rm supp}~a_k \subset 2^{k }Q_{1} ,$ and
$$  2^{n  b }  \left\{\sum
\limits_{k=0}^{\infty}2^{-\bar{p}{ank}}\right\}^{1/\bar{p}}=C<\infty$$
since $a>0$. That is
$\|M \|^{ }_{_{ }BL^{p,   \alpha}}\leq C .$
 Thus,    Theorem 2 is proved.
\\

\par {\bf Proposition 1 }   Let
 $   0< p < \infty$, $  -n < \alpha < n(p-1)$
 .
  Then,
$$\|M _{ {\mathcal{C}} _0 } h\|_{ BL^{p,   \alpha}}\leq
C,  $$
for a $ (  p, \alpha)-$block $h$
,
where constant $C$ is independent of $h
$
.
\par {\bf Proof  }
By the molecular theorem, it  suffices to check that $M _{ {\mathcal{C}} _0 }h$ is a  $ (p,  \alpha, \varepsilon
)$-molecular centered at $x_0$
 for every $_{ }(  p, \alpha)$- block $h$  centered at $x_0$ and $\Re (M _{ {\mathcal{C}} _0 }h)\leq C $, where $C$ is independent of $h$.
 Given a $_{ }( p, \alpha)-$block
$h$ with
 supp $h\subset  Q_{ 0}  $ in $ {\mathcal{C}}    $$_{1}$ centered at $x_0$. We see that $\|M _{ {\mathcal{C}} _0 }h\|_{L^{\infty}}\leq \|h\|_{L^{\infty}} \leq
|Q_{0} |^{ -1/p-\alpha/np}$.
    Since  $1
 -1/p-\alpha/np>0$, we can choose $\varepsilon$ such that
$$0<\varepsilon <\min\{1 -1/p-\alpha/np,1 \},$$
 and set $$ a=1-1/p-\alpha/np-\varepsilon  ~~{\rm  and }~~ b=1 -\varepsilon. $$ Clearly, $a>0,b>0$ and $b-a= (  n+\alpha)/np> 0$.
To get $\Re (M _{ {\mathcal{C}} _0 }h)\leq C $,
 it  suffices to prove   that
$$J=:
\|M _{ {\mathcal{C}} _0 }h(x)|x -x_0 | ^{nb}\|^{}_{L^{\infty} }
\leq C |Q_{0} |^b\|h \|^{}_{L^{\infty}   }.\eqno (3)$$
Write
\begin{eqnarray*}
J\leq\|M _{ {\mathcal{C}} _0 }h(x)|x -x _0| ^{nb}\chi_{2Q_{0}}(x)\|^{}_{L^{\infty}
}
 +\|M _{ {\mathcal{C}} _0 }h(x)|x -x _0| ^{nb}\chi_{(2Q_{0})^c}(x)\|^{}_{L^{\infty}
}
  =:  J_1+J_2.
\end{eqnarray*}
 For $J_1$,
 we  have   that
\begin{eqnarray*}
J_1
 \leq   \|M _{ {\mathcal{C}} _0 }h
\|^{}_{L^{\infty}   }\| |x -x _0| ^{nb}\chi_{2Q_{0}}(x)\|^{}_{L^{\infty}
}
  \leq  C |Q_{ 0}|^{b} \|h \|^{}_{L^{\infty}   }
\end{eqnarray*}
 since $b>0$ and   $L^{\infty}$ boundedness of $M _{ {\mathcal{C}} _0 }$.
For  $J_2,$    we see that
$$
M _{ {\mathcal{C}} _0 }h(x)=\sup_{x\in C
 \in {{\mathcal{C}} _0 }, {\rm diam }(C)<1 } \frac{1}{|C|}\int_C |h(y)|dy=\sup_{x\in C
 \in {{\mathcal{C}} _0 }, {\rm diam} (C)<1 } \frac{1}{|C|}\int_{C\cap Q_0} |h(y)|dy=0
$$
for $x \in (2Q_{0})^c $, since $C\cap Q_0=\emptyset $. So, $J_2=0$. Thus, (3) holds.    Proposition 1 follows.

\par   Theorem 1 follows easily from Proposition 1.

\par
  Let $0<p<\infty$ and   $-n < \alpha < n(p-1)$,
  it is known that
        $_{ } {B}L^{p,   \alpha} \subset L^{\frac{np}{n+\alpha}},$
   see \cite{L}.   Then,

\par {\bf  Corollary 1} Let $0<p<\infty$ and $  - n < \alpha < n(p-1)  $.  $M _{ {\mathcal{C}} _0 } $ is bounded from
$
BL^{p,   \alpha} $ to  $L^{\frac{np}{n+\alpha}}$.


 {\bf
 3. The Differentiation Theorems}   By Corollary 1 and a standard arguments, we have that

\par {\bf Theorem 3} Let $0<p<\infty$ and $  - n < \alpha < n(p-1)  $.  Let $B$ be a subspace of  $BL^{p,   \alpha}$. If (1) holds as $   {\mathcal{C}} =    {\mathcal{C}} _0     $ for the functions in a dense subspace $D$  of
$
B$, then (1) holds for all functions in
$
 B $.

\par $BC^{p,   \alpha}$ and $B^{p,   \alpha}$ are in $BL^{p,   \alpha}$. At the same time, we see that

\par {\bf Property 1} Let   $0<p<\infty, $ and $ -n<\alpha<\infty.$   $C^{ }_c$
 is
dense in $BC^{p,   \alpha}.$

\par{\bf Proof} Let $f\in BC^{p,   \alpha}$, i.e. $f(x)=\sum_{k=1}^{\infty} \lambda_ka_k(x)
$ where each $a_k$ is a $(p, \alpha)$-continuous block with supp $
a_k\subseteq Q_k$ and $ \sum_{k=1}^{\infty}
 |\lambda_k|^{\bar{p}}<\infty $. Then for any $\varepsilon>0,$ there exists
 an
 $i_0$ such that
$
 \sum_{k=i_0+1}^{\infty}
 |\lambda_k|^{\bar{p}}<\varepsilon^{\bar{p}}.
 $
Let $f_{i_0}(x)=\sum_{k=1}^{i_0} \lambda_ka_k(x) $, we see that $f_{i_0} \in C^{ }_c$,
and
\begin{eqnarray*}
\|f-f _{i_0} \|_{BC^{p,   \alpha}} ^{\bar{p}} =
\|\sum_{k=i_0+1}^{\infty}
\lambda_ka_k \|_{BC^{p,   \alpha}} ^{\bar{p}}
   \leq
 \sum_{k=i_0+1}^{\infty}
 |\lambda_k|^{\bar{p}}
 <\varepsilon ^{\bar{p}}.
 \end{eqnarray*}
 Thus, Property 1 holds.
 \par Let    $    {\mathcal{C}} _{\rm R}   $  be the collection of   arbitrary rectangles in  ${\bf R}^n  $. it is known that (1) holds as $   {\mathcal{C}} =    {\mathcal{C}} _{\rm R}    $  for all functions in
$C_c$. So, we have that

 \par {\bf  Corollary 2} Let $0<p<\infty$ and $  - n < \alpha < n(p-1)  $.       (1) holds as $   {\mathcal{C}} =    {\mathcal{C}} _{\rm R}     $  for all functions in
$BC^{p,   \alpha}$.

 \par Let $FB  ^{p,   \alpha} $   be the functions class consists of finite   linear combinations of $(p, \alpha) $- characteristic blocks. It is easy to see that $FB^{p  ,\alpha}$ is dense in $B^{p,   \alpha}$. While (1) holds as $   {\mathcal{C}} =    {\mathcal{C}} _0     $  for all functions in
$FB^{p  ,\alpha}$, in fact, for  each $(p,  \alpha) $-characteristic block $a_k$ with support set $Q$ (cube), the equation in (1) holds for   $x\in Q\backslash \partial Q$ (inner point of $Q$). So,

  \par {\bf  Corollary 3} Let $0<p<\infty$ and $  - n < \alpha < n(p-1)  $.        (1) holds as $   {\mathcal{C}} =    {\mathcal{C}} _0     $  for all functions in
$B^{p,   \alpha}$.

 {\bf  4.   A count example  }
  Denote {\bf Q} the set of rational numbers.
  Denote  $I_i$ the  closed interval $ [i,i+1], i\in {\bf Z}  $. Let ${\bf Q}_{i,j}=\{
(r,q)\in I_i\times I_j:r,q\in {\bf Q} \}$.  Let $E_{r,q}^{i,j}$ be the closed cube with sides parallel to the axes, with side length 1 and with a vertex $(r,q)$ in the upper right corner.
   Let
 $$f_{i,j}(x)=\sum_{(r,q)\in {\bf Q}_{i,j}}\lambda_{r,q}^{i,j}
  \chi_{E_{r,q}^{i,j}}(x),$$
with each $ \lambda_{r,q}^{i,j} >0$ for even $j $ and all $i\in {\bf Z}$ and $ \lambda_{r,q}^{i,j} <0$ for odd  $j $ and all $i\in {\bf Z}$, and $ \sum_{(r,q)\in {\bf Q}_{i,j}}|\lambda_{r,q}^{i,j}|^{\bar{p}} <\infty.$ Let
$$g(x)=\sum_{i,j} f_{i,j}(x).$$

\par It is easy to see that $g(x) \in B^{p,   \alpha}.$ So (1) holds for $g$ as $   {\mathcal{C}} =    {\mathcal{C}} _0     $.

\par {\bf Proposition 2} There exists
$g(x)  $ which is   continuous nowhere on ${\bf R}^2$ such that for   which (1) holds as $   {\mathcal{C}} =    {\mathcal{C}} _0     $.

\par {\bf Proof} Need only to prove that $g(x) $ above  is   continuous nowhere on ${\bf R}^2$. Taking $(r_i,q_j) \in {\bf Q}_{i,j}$ for each pair $(i,j)$, we see that
$$(r_i,q_j)\in \cap _{ q_j \leq q \leq   j+1 ,  q\in {\bf Q}  }E_{r_i,q}^{i,j},
~~(r_i,q_j)\in \cap _{ q_j \leq q \leq   j+1 ,  q\in {\bf Q}  }E_{r_i+1,q}^{i+1,j},$$
$$(r_i,q_j)\in \cap _{j+1 \leq q \leq   q_j+1 ,  q\in {\bf Q}  }E_{r_i ,q}^{i ,j+1},~~(r_i,q_j)\in \cap _{j+1 \leq q \leq   q_j+1 ,  q\in {\bf Q}  }E_{r_i+1,q}^{i+1,j+1},
$$
it follows that
$$f_{i,j}((r_i,q_j))=\sum_{ q_j \leq q \leq   j+1 ,  q\in {\bf Q}  }\lambda_{r_i,q}^{i,j},$$
$$f_{i+1,j}((r_i,q_j))=\sum_{ q_j \leq q \leq   j+1 ,  q\in {\bf Q}  }\lambda_{r_i+1,q}^{i+1,j},$$
$$f_{i ,j+1}((r_i,q_j))=\sum_{ j+1 \leq q \leq   q_j+1 ,  q\in {\bf Q}  }\lambda_{r_i ,q}^{i ,j+1},$$
$$f_{i+1,j+1}((r_i,q_j))=\sum_{ j+1 \leq q \leq   q_j+1 ,  q\in {\bf Q}  }\lambda_{r_i+1,q}^{i+1,j+1}.$$
Noticing
$$ g((r_i,q_j)) =f_{i,j}((r_i,q_j))+f_{i+1,j}((r_i,q_j)) +f_{i ,j+1}((r_i,q_j)) +f_{i+1,j+1}((r_i,q_j)), $$
we have for $\tilde{q}_j>q_j$ that

$
g((r_i,\tilde{q}_j))-g((r_i,q_j))=$
$$-\sum_{ q_j \leq q <  \tilde{ q}_j ,  q\in {\bf Q}  }\lambda_{r_i,q}^{i,j}
-\sum_{ q_j \leq q <  \tilde{ q}_j ,  q\in {\bf Q}  }\lambda_{r_i+1,q}^{i+1,j}
+\sum_{ q_j +1\leq q <  \tilde{ q}_j +1,  q\in {\bf Q}  }\lambda_{r_i ,q}^{i ,j+1}
+\sum_{ q_j +1\leq q <  \tilde{ q}_j +1,  q\in {\bf Q}  }\lambda_{r_i+1,q}^{i+1,j+1}
,
$$
it follows that
$$
|g((r_i,\tilde{q}_j))-g((r_i,q_j))| > |\lambda_{r_i,q_j}^{i,j}|.$$
It follows that  $f(x)$ is not continuous at $(r_i,q_j)$. Proposition 2 follows.

{\bf  5.   Two remarks  }

 {\bf Remark 1}  Let $    {\mathcal{C}}_l    $$ =\{ cubes  ~ Q
    :|Q|\geq l^n>0   \}$. Let $ BL^{p,   \alpha}_l,  BC^{p,   \alpha}_l$ and $  B^{p,   \alpha}_l$ be the classes
 replacing $    {\mathcal{C}}_1    $ in Definition 1 and 2 by $    {\mathcal{C}}_l    $ . By a simple dilation transform, we see that Theorem 3, Corollary 2 and 3, and Proposition 1 hold  for $l>0 $ if $ BL^{p,   \alpha} ,  BC^{p,   \alpha} $ and $  B^{p,   \alpha} $ are replaced by $ BL^{p,   \alpha}_l,  BC^{p,   \alpha}_l$ and $  B^{p,   \alpha}_l$, respectively.  Let
 $
M _{{ {\mathcal{C}} _0 },l}f(x)=\sup_{x\in C \in { {\mathcal{C}} _0 }, {\rm diam}(C)<l } \frac{1}{|C|}\int_C |f(y)|dy,
$ it is easy to see from the proofs that Theorem 1 and 2, Proposition 1 and Corollary 1  hold  for $
M _{{\mathcal{C} }_0, l}$ and $ BL^{p,   \alpha}_l$ when  $l>0 $.

 {\bf Remark 2}  Let $ BL^{p,   \alpha}_{{\rm R}},  BC^{p,   \alpha}_{{\rm R}}$ and $  B^{p,   \alpha}_{{\rm R}}$ be the classes
 replacing $    {\mathcal{C}}_1    $ in Definition 1 and 2 by $    {\mathcal{R}}_0    $$ =\{ rectangles  ~ R ~with ~sides ~length \geq 1 ~and~   fixed~eccentricity
          \}$, (the eccentricity of a rectangle is the ratio of its longest size over its shortest size). Then, Theorem 1-3, Corollary 1-3 and Proposition 1 hold still if $ BL^{p,   \alpha} ,  BC^{p,   \alpha} $ and $  B^{p,   \alpha} $ are replaced by $ BL^{p,   \alpha}_{{\rm R}},  BC^{p,   \alpha}_{{\rm R}}$ and $  B^{p,   \alpha}_{{\rm R}}$, respectively. This can be seen by the  facts:
 $$R\subset Q {\rm~~ and~~ } |R|^{-\alpha/pn-1/p  }\leq (d_{\rm min}^n)^{-\alpha/pn-1/p  } =c^{n(\alpha/pn+1/p)}|Q|^{-\alpha/pn-1/p  } , $$
 since $-\alpha/pn-1/p <0$, where $R\in {\mathcal{R}}_0 $  with    the shortest side length $d_{\rm min} $  and   eccentricity $c$, $Q$ is the cube with sides length $cd_{\rm min} $.

\par {\bf Acknowledgements}  This work was supported in part by NSF of China grant  11171280

{ }
\par
Shunchao Long
\par Department of Mathematics,
\par Xiangtan University,
   Hunan, 411105  China
 \par E-mail: sclong@xtu.edu.cn

\end{document}